\documentclass[12pt,reqno]{amsart}
\usepackage{amsmath,amsthm,amsfonts,amssymb}
\usepackage[ansinew]{inputenc}
\usepackage{mathrsfs}

\newcommand{\R}{\mathbf{R}}

\newcommand{\Z}{\mathbf{Z}}

\newcommand{\fg}{{\mathfrak g}}
\newcommand{\fh}{{\mathfrak h}}
\newcommand{\fk}{{\mathfrak k}}
\newcommand{\fp}{{\mathfrak p}}
\newcommand{\fm}{{\mathfrak m}}
\newcommand{\ft}{{\mathfrak t}}
\newcommand{\fa}{{\mathfrak a}}
\newcommand{\fz}{{\mathfrak z}}

\newcommand{\SO}{{\rm SO}}
\newcommand{\U}{{\rm U}}
\newcommand{\rmS}{{\rm S}}
\newcommand{\rmO}{{\rm O}}

\newcommand{\pr}{{\rm pr}}

\newcommand{\fso}{{\mathfrak{so}}}
\newcommand{\fu}{{\mathfrak{u}}}

\newcommand{\Exp}{{\rm Exp}}
\newcommand{\ad}{{\rm ad}}
\newcommand{\Ad}{{\rm Ad}}
\newcommand{\Vol}{{\rm Vol}}

\newtheorem{Thm}{Theorem}[section]

\newtheorem{Prop}[Thm]{Proposition}
\newtheorem{Cor}[Thm]{Corollary}
\newtheorem{Lemma}[Thm]{Lemma}

\newtheorem*{Thm2}{Theorem}

\theoremstyle{remark}
\newtheorem*{Rem}{Remark}

\begin{document}

\title[Geometry of Hermann Actions]{On the Geometry of the Orbits of
Hermann Actions}

\author[Goertsches]{Oliver Goertsches}
\address[Goertsches]{%
    Max-Planck-Institut f\"ur Mathematik,
    Vivatsgasse 7,
   53111 Bonn,
   Germany
} \email{olig@mpim-bonn.mpg.de}

\author[Thorbergsson]{Gudlaugur Thorbergsson}
\address[Thorbergsson]{%
Mathematisches Institut, Universit\"at zu K\"oln, Weyertal 86 -
90, 50931 K\"oln,
Germany %
} \email{gthorbergsson@mi.uni-koeln.de}


\begin{abstract}
We investigate the submanifold geometry of the orbits of Hermann
actions on Riemannian symmetric spaces. After proving that the
curvature and shape operators of these orbits commute, we
calculate the eigenvalues of the shape operators in terms of the
restricted roots. As applications, we get a formula for the
volumes of the orbits and a new proof of a Weyl-type integration
formula for Hermann actions.
\end{abstract}

\thanks{The second author was supported by the DFG-Schwerpunkt {\it
Globale Differential\-geometrie}.}

\maketitle

\tableofcontents

\section{Introduction and Results}
An isometric action of a compact Lie group on a Riemannian
manifold $M$ is called {\it polar} if it admits a {\it section},
i.e.~a connected submanifold $\Sigma$ of $M$ that meets all
orbits perpendicularly at each point of intersection. If the
section is flat, the action is called {\it hyperpolar}.

In this paper, $M=G/K$ will denote a Riemannian symmetric space
of compact type. As the classification of hyperpolar actions on
irreducible symmetric spaces of compact type \cite{Kollross}
shows, all examples of such actions of cohomogeneity at least two
are orbit equivalent to the so-called {\it Hermann actions},
i.e.~actions of symmetric subgroups of $G$. Recall that a
subgroup $H\subset G$ is called {\it symmetric} if there exists an
involutive automorphism $\sigma:\fg\to \fg$ with fixed point
algebra $\fh$.

In Section \ref{sec_proof}, we prove the following theorem.
\begin{Thm2}  Let $H\subset G$ be a symmetric subgroup, $p\in M$ regular
and $v,w\in \nu_pHp$. Then the tangent space $T_pHp$ is an
invariant subspace of the curvature operator $R_v(x)=R(x,v)v$ and
the restriction of $R_v$ to $T_pHp$ commutes with the shape
operator $A_w$ of $Hp$.
\end{Thm2}

Therefore the curvature and shape operators of $Hp$ can be
simultaneously diagonalized. The eigenspaces of the curvature
operators are given by the root spaces of $M$; more precisely, a
coarser version of the root space decomposition obtained by
regarding only the restrictions of the roots to the tangent space
of the section is relevant here -- see Section
\ref{sec_RootSpaceDecomp}. As a corollary of the above theorem,
we obtain that for singular orbits, the restricted curvature
operator $R_v$ commutes with the shape operator $A_w$ if $v$ and
$w$ lie in the same section.

In Section \ref{sec_EV_ComInv}, we restrict ourselves to the case
where $H$ can be conjugated in such a way that the involutions
corresponding to $H$ and $K$ commute\footnote{Note that in this
case the triple $(G,H,K)$ is called a {\it symmetric triad} in
\cite{ConlonTop}.}, which is possible except in a few cases
\cite{Conlon}. We can completely determine the eigenspaces of the
shape operators in terms of the restrictions of the roots
(Theorem \ref{thm_ShapeOp}), thereby generalizing \cite{Verhoczki}
where the case $H=K$ is treated.

In the general case, which is treated in Section
\ref{sec_EV_General}, we can show how the eigenvalues of $A_v$
change if the normal direction $v$ is varied (Proposition
\ref{prop_ShapeOpGeneral}).

Using the methods of \cite{CNV}, where the case $H=K$ is treated,
we calculate in Section \ref{sec_Appl} the volumes of the
principal orbits; furthermore we reprove a Weyl-type integration
formula for actions of Hermann type (\cite{F-J}, which is a
generalization of Theorem I.5.10 of \cite{Helgason_GGA}) using
our calculations of the shape operators.

We would like to remark that, with slight modifications, our
results are also true in the noncompact case, the only difference
being some sign changes and some replacements of trigonometric
functions by hyperbolic ones. Nevertheless, for better
readability, we will present the proofs only for the compact
case. Note that in the noncompact case, $H$ can always be
conjugated in such a way that the two involutions commute, see
\cite{Berger}, Lemma 10.2. The shape operators in the noncompact
case are also calculated in \cite{Koike}, but in a completely
different way.

\section{Preliminaries}
Let $M=G/K$ be a Riemannian symmetric space of compact type and
set $p=eK$. Then $G$ is a semisimple compact Lie group, and we
assume the metric on $M$ to be induced by the Killing form of
$G$. The Lie algebra $\fg$ can be identified with the Lie algebra
of Killing vector fields on $M$, with the bracket being the
negative of the bracket of the Killing vector fields. Considering
the Cartan decomposition $$\fg=\fk\oplus \fm$$ and using the
identification of $\fg$ with the Killing vector fields, we have
\begin{equation}\fk=\{X\in \fg\mid X(p)=0\}\quad\text{and}\quad \fm=\{X\in
\fg\mid (\nabla X)_{p}=0\},\label{eqn_G_Killing}
\end{equation}see \cite{Sakai}, Lemma 6.8. The Killing vector fields in
$\fm$ are those induced by transvections along geodesics through
$p$.

If $X,Y,Z\in \fm$, we can express the curvature of $M$ at the
point $p$ by
\begin{equation}
R(X(p),Y(p))Z(p)=-[[X,Y],Z](p). \label{eqn_CurvSymmSpace}
\end{equation}
Note that this equality remains valid if we assume only two of
the Killing vector fields to be induced by transvections -- if
e.g.~$X\in \fg$ is arbitrary, this follows from
$[[\fk,\fm],\fm]\subset [\fm,\fm]\subset \fk$.

Let now $\fa\subset \fm$ be a maximal abelian subalgebra, denote
the set of restricted roots by $\Delta$ and a choice of positive
roots by $\Delta^+$. Then the corresponding root space
decomposition of $M$ is
\begin{equation}
\label{eqn_RootSD_aIntr} \fk=\fz_\fk(\fa)\oplus \sum_{\alpha\in
\Delta^+} \fk_\alpha\quad \text{ and }\quad \fm=\fa\oplus
\sum_{\alpha\in \Delta^+} \fm_\alpha,
\end{equation}
where
\begin{equation}
\fk_\alpha=\{X\in \fk\mid \ad_W^2(X)=-\alpha(W)^2X \text{ for all
}W\in \fa\} \label{eqn_K_Alpha}
\end{equation}
and
\begin{equation}
\fm_\alpha=\{X\in \fm\mid \ad_W^2(X)=-\alpha(W)^2X \text{ for all
}W\in \fa\}.\label{eqn_M_Alpha}
\end{equation}
We call $X\in \fk_\alpha$ and $\fm_\alpha$ {\it related} if
$[W,X]=-\alpha(W)Y$ and $[W,Y]=\alpha(W)X$ for all $W\in \fa$ (see
\cite{Loos}, p.~61). For any $X\in \fk_\alpha$ there exists a
related vector $Y\in \fm_\alpha$, and vice versa; in particular,
the vector spaces
$\fk_\alpha$ and $\fm_\alpha$ are isomorphic.\\

For $v\in T_pM$, the {\it curvature operator} $R_v$ is defined to
be the endomorphism of $T_pM$ given by $R_v(u)=R(u,v)v$.

The {\it shape operator} $A_\xi:T_pN\to T_pN$ of a submanifold
$N\subset M$ in the normal direction $\xi\in \nu_pN$ is defined
as $A_\xi x=-(\nabla_x \xi)^T$; with this choice of sign, a
Jacobi field $J$ along the normal geodesic $\gamma$ in direction
$\xi$ is an $N$-Jacobi field if and only if $J(0)\in T_pN$ and
$J'(0)+A_\xi J(0)\in \nu_pN$.

Let $H\subset G$ act on $M=G/K$. If $p\in M$ is regular, the fact
that the slice representation at $p$ is trivial implies that we
can extend normal vectors to well-defined $H$-equivariant normal
vector fields on $Hp$. If the $H$-action is polar, these are
automatically parallel with respect to the normal connection, see
\cite{PalaisTerng}, Theorem 5.6.7. For any such vector field
$\xi$, we thus get
\begin{equation}\label{Eqn_ShapeOpEquivNVF}
A_{\xi(p)}x=-\nabla_x\xi.
\end{equation}
The exponential map of $G$ will be denoted by $\exp$ and the one
of $M$ by $\Exp$. For parallel translation along a curve
$\gamma$, we will write
${}_\gamma||_{t_0}^{t_1}:T_{\gamma(t_0)}M\to T_{\gamma(t_1)}M$.

\section{Proof of the Main Theorem}\label{sec_proof}
In this section, we will first prove the following theorem about
the principal orbits of Hermann actions; the singular orbits will
be dealt with in Corollary \ref{Cor_CommuteSingular}.

\begin{Thm}\label{Thm_Commute} Let $H\subset G$ be a symmetric subgroup,
$p\in M$ regular and $v,w\in \nu_pHp$. Then the tangent space
$T_pHp$ is an invariant subspace of the curvature operator
$R_v(x)=R(x,v)v$ and the restriction of $R_v$ to $T_pHp$ commutes
with the shape operator $A_w$ of $Hp$.
\end{Thm}
First we need a lemma.
\begin{Lemma}\label{Lem_H_Xi_Commute}
Let $X\in \fh$ and $\xi$ be an $H$-equivariant normal vector
field on $Hp$. Then $[X,\xi]=0$.
\end{Lemma}
\begin{proof} Let $\gamma(t)=\Exp(t\xi(p))$. Then we have
\begin{align*}
\nabla_{X(p)}\xi&=\left.\frac{\nabla}{ds}\right|_{s=0}\xi(\exp(sX)\cdot
p)=\left.\frac{\nabla}{ds}\right|_{s=0}\left.\frac{d}{dt}\right|_{t=0}\exp(sX)\cdot\gamma(t)\\
&=\left.\frac{\nabla}{dt}\right|_{t=0}\left.\frac{d}{ds}\right|_{s=0}\exp(sX)\cdot
\gamma(t)=\left.\frac{\nabla}{dt}\right|_{t=0}X(\gamma(t))=\nabla_{\xi(p)}X.
\end{align*}
\end{proof}

\begin{proof}[Proof of Theorem \ref{Thm_Commute}] The invariance of
$T_pHp$ under $R_v$ follows from the fact that for any $x\in
T_pHp$ and any $u\in\nu_pHp$, we have
$$\left<R(x,v)v,u\right>=-\left<R(v,u)v,x\right>=0$$ because the action is
polar.

Without loss of generality, we may assume that $p=eK$. Note that
then $\ft:=\fh^\bot\cap \fm\cong \nu_pHp$ is an abelian
subalgebra of $\fg$.

Let $\vartheta\in \ft$ be the Killing vector field on $M$ with
$\vartheta(p)=v$ and $(\nabla \vartheta)(p)=0$ and let $\xi$ be
the $H$-equivariant normal vector field on $Hp$ with $\xi(p)=w$.
Furthermore, let $\xi'\in \ft$ be the Killing vector field with
$\xi'(p)=w$ and $(\nabla\xi')(p)=0$ and set $g(t)=\exp(t\xi(p))$.
For $X\in \fh$, we have
\begin{align*}
R_v(A_w&(X(p)))\overset{(\ref{Eqn_ShapeOpEquivNVF})}{=}-R(\nabla_{X(p)}
\xi,v)v\overset{\ref{Lem_H_Xi_Commute}}{=}-R(\nabla_{\xi(p)} X,v)v\\
&=-\nabla_{\xi(p)}(R(X,\vartheta)\vartheta)+R(X(p),\underbrace{\nabla_{\xi(p)}\vartheta}_{=0})v+R(X(p),v)\underbrace{\nabla_{\xi(p)}
\vartheta}_{=0}\\
&\!\!\!\!\!\!\!\!\!\!\!
\text{(}\nabla R=0\text{ on a symmetric space)}\\
&=-\left.\frac{\nabla}{dt}\right|_{t=0}
R(X(g(t)p),\vartheta(g(t)p))\vartheta(g(t)p)\\
&=-\left.\frac{\nabla}{dt}\right|_{t=0}d(g(t))(R(\Ad_{g(t)^{-1}}X,\Ad_{g(t)^{-1}}\vartheta)\Ad_{g(t)^{-1}}\vartheta)\\
&=-\left.\frac{\nabla}{dt}\right|_{t=0}d(g(t))
R(\Ad_{g(t)^{-1}}X,\vartheta(p))\vartheta(p)\\
&\!\!\!\!\!\!\!\!\!\!\!
\text{(}\ft\text{ is abelian)}\\
&=-\left.\frac{\nabla}{dt}\right|_{t=0}{}_\gamma||_0^t
R(\Ad_{g(t)^{-1}}X,\vartheta(p))\vartheta(p)\\
&=-\left.\frac{d}{dt}\right|_{t=0}R(\Ad_{g(t)^{-1}}X,\vartheta(p))\vartheta(p)\\
&=\left.\frac{d}{dt}\right|_{t=0}[[\Ad_{g(t)^{-1}}X,\vartheta],\vartheta](p)\\
&\!\!\!\!\!\!\!\!\!\!\! \text{(}\vartheta\text{ is induced by
transvections; see
(\ref{eqn_CurvSymmSpace}))}\\
&= [[[X,\xi'],\vartheta],\vartheta](p)\\
&=[[[X,\vartheta],\xi'],\vartheta](p)+[[X,\underbrace{[\xi',\vartheta]}_{=0}],\vartheta](p)\\
 &=[[[X,\vartheta],\vartheta],\xi'](p)\\
&=-\underbrace{\nabla_{[[X,\vartheta],\vartheta](p)}\xi'}_{=0\,
(\xi'\in \fp)}+\nabla_{\xi(p)}[[X,\vartheta],\vartheta]\\
&\!\!\!\!\!\!\!\!\!\!\!
\text{(note the sign of the bracket in }\fg \text{)}\\
&=-\nabla_{R(X(p),v)v}\xi+[[[X,\vartheta],\vartheta],\xi](p)\qquad\,\,\quad\qquad\qquad\qquad\\
&=A_w(R_v(X(p)))+[[[X,\vartheta],\vartheta],\xi](p).
\end{align*}
Since $X$ and $\vartheta$ are Killing vector fields on $M$, we
have $[[X,\vartheta],\vartheta]\in
[[\fh,\fh^\bot],\fh^\bot]\subset [\fh^\bot,\fh^\bot]\subset \fh$,
where for the last inclusion we used that $H$ is a symmetric
subgroup of $G$. Hence, Lemma \ref{Lem_H_Xi_Commute} implies that
$[[[X,\vartheta],\vartheta],\xi](p)=0$.
\end{proof}

For the singular orbits, we have the following corollary.
\begin{Cor}\label{Cor_CommuteSingular} Let $H\subset G$ be a symmetric
subgroup and $p\in M$ arbitrary. Then for all $v\in \nu_pHp$, the
tangent space $T_pHp$ is invariant under the curvature operator
$R_v$. If $\Sigma$ is a section of the $H$-action passing through
$p$ and $v,w\in T_p\Sigma$, then the restriction of $R_v$ to
$T_pHp$ commutes with the shape operator $A_w$ of $Hp$.
\end{Cor}
\begin{proof} Conjugate $H$ such that $p$ is the origin. Let $\vartheta\in
\fh^\bot\cap \fm$ be the Killing vector field induced by
transvections with $\vartheta(p)=v$. Then for any $X\in \fh$,
equation (\ref{eqn_CurvSymmSpace}) yields
$$R(X(p),\vartheta(p))\vartheta(p)=-[[X,\vartheta],\vartheta](p)\in
T_pHp$$ because $\fh$ is a symmetric subgroup.

Then the operators commute because the regular points in $\Sigma$
are dense in $\Sigma$.
\end{proof}

\begin{Cor}\label{Cor_AllCommute} Let $H\subset G$ be a symmetric
subgroup. Then $\{A_v,R_v\mid v\in T_p\Sigma\}$ is a commuting
family of endomorphisms of $T_pHp$, where $\Sigma$ is a section
passing through $p$.
\end{Cor}
\begin{proof} We first assume that $p$ is regular. Then the $R_v$ commute
because $\nu_pHp$ is abelian as one sees by combining
(\ref{eqn_CurvSymmSpace}) with the Jacobi identity. The Ricci
equation implies
\begin{align*}
\left<[A_v,A_w]x,y\right>&=\left<R(x,y)v,w\right>-\left<R^\bot(x,y)v,w\right>\\
&=\left<R(v,w)x,y\right>=0
\end{align*}
for all $x,y\in T_pHp$ and all $v,w\in \nu_pHp$ since the normal
bundle of $Hp$ is flat. Therefore, the shape operators of $Hp$
commute.

If $p$ is not regular, the claim then follows from the fact that
the regular points in $\Sigma$ are dense in $\Sigma$.
\end{proof}

\begin{Rem}
A submanifold $N$ of a Riemannian manifold $M$ is called {\it
curvature-adapted} if $T_pN$ is invariant under the curvature
operator $R_u$ and if the restriction of $R_u$ to $T_pN$ commutes
with the shape operator $A_u$ of $N$ for any $p\in N$ and all
$u\in \nu_p M$. Theorem \ref{Thm_Commute} and Corollary
\ref{Cor_CommuteSingular} immediately imply that all the orbits
of the $H$-action are curvature-adapted submanifolds of $M$.
\end{Rem}

\section{Adapted Root Space Decomposition}
\label{sec_RootSpaceDecomp}

The root space decomposition adapted to the $H$-action has
already been described in \cite{F-J}; we include it for the
convenience of the reader.

Let $H\subset G$ act on $M$ in a hyperpolar fashion; later we
will assume that it is a symmetric subgroup such that the
corresponding involution commutes with the one of $\fk$. Consider
the Cartan decomposition $$\fg=\fk\oplus \fm$$ and choose a
maximal abelian subalgebra $\ft\subset \fh^\bot \cap \fm$. Then
$\Sigma=\Exp(\ft)\subset M$ is a section for the $H$-action on
$M$. Note that $\Sigma$ is a torus since we assume the metric to
be induced by the Killing form on $\fg$; see \cite{HPTT}, Theorem
2.3. Note that $$\fm=\pr_\fm\fh\oplus (\fh^\bot\cap \fm)\cong
T_pHp\oplus \nu_pHp.$$ Let further $\fa$ be a maximal abelian
subalgebra of $\fm$ containing $\ft$. The set of restricted roots
of $M$ with respect to $\fa$ shall be denoted by $\Delta$ and a
choice of positive roots by $\Delta^+$.

Consider the corresponding root space decomposition
\begin{equation}
\label{eqn_RootSD_a} \fk=\fz_\fk(\fa)\oplus \sum_{\alpha\in
\Delta^+} \fk_\alpha\quad \text{ and }\quad \fm=\fa\oplus
\sum_{\alpha\in \Delta^+} \fm_\alpha,
\end{equation}
where the root spaces $\fk_\alpha$ and $\fm_\alpha$ are given by
(\ref{eqn_K_Alpha}) and (\ref{eqn_M_Alpha}).

In the following, the restrictions of the roots to $\ft$ will be
of greater importance than the roots themselves. We define
$$\Delta_\ft=\{\alpha|_\ft\mid \alpha\in \Delta\}\setminus\{0\}.$$ Note
that for two roots $\alpha,\alpha'\in \Delta^+$, it is possible
that $\alpha|_\ft=-\alpha'|_\ft$ (see the example in Section
\ref{sec_Example}). Therefore, we let $\Delta_\ft^+\subset
\Delta_\ft$ be the set of nonzero restrictions of elements in
$\Delta^+$ to $\ft$, but if this occurs, we include only one of
$\alpha|_\ft$ and $\alpha'|_\ft$ in $\Delta_\ft^+$. For any
$\beta\in \Delta^+_\ft$, we set
$$\fk_\beta^\ft=\{X\in \fk\mid \ad_w^2(X)=-\beta(w)^2X \text{ for all }
w\in \ft\}$$ and
$$\fm_\beta^\ft=\{X\in \fm\mid \ad_w^2(X)=-\beta(w)^2X \text{ for all }
w\in \ft\}.$$

\begin{Lemma}
\label{lem_RelationRootSDs} For any $\beta\in \Delta^+_\ft$, we
have
$$\fk_\beta^\ft=\sum_{\alpha\in \Delta^+:\, \alpha|_\ft=\pm\beta}
\fk_\alpha\quad \text{ and }\quad  \fm_\beta^\ft=\sum_{\alpha\in
\Delta^+:\, \alpha|_\ft=\pm\beta} \fm_\alpha.$$ Furthermore,
$$\fz_\fm(\ft)=\fa\oplus \sum_{\alpha\in \Delta^+:\, \alpha|_\ft=0}
\fm_\alpha \quad \text{ and } \quad
\fz_\fk(\ft)=\fz_\fk(\fa)\oplus \sum_{\alpha\in \Delta^+:\,
\alpha|_\ft=0} \fk_\alpha.$$
\end{Lemma}
\begin{proof}
We only prove the first equality. For all $\beta\in \Delta^+_\ft$
and all $\alpha\in \Delta^+$ such that $\alpha|_\ft=\pm\beta$, we
have $\fk_\alpha\subset \fk_\beta^\ft$ by definition. The desired
equality then follows from (\ref{eqn_RootSD_a}).
\end{proof}

Consequently we have the decompositions
\begin{equation}
\label{eqn_RootSD_t} \fk=\fz_\fk(\ft)\oplus \sum_{\beta\in
\Delta^+_\ft}\fk_\beta^\ft \quad \text{ and }\quad
\fm=\fz_\fm(\ft)\oplus\sum_{\beta\in \Delta^+_\ft}\fm_\beta^\ft.
\end{equation}

From now on until the end of Section \ref{sec_EV_ComInv}, we
assume that $H$ and $K$ are symmetric subgroups of $G$
corresponding to commuting involutions. Then the two Cartan
decompositions
\begin{equation}
\fg=\fk\oplus \fm =\fh\oplus\fp
\end{equation}
are compatible in the sense that
\begin{equation}
\label{eqn_CartanRefined} \fg=\fk\cap \fp\oplus \fk\cap \fh\oplus
\fm\cap \fh \oplus \fm\cap \fp.
\end{equation}

\begin{Lemma}
For all $\beta\in \Delta_\ft^+$, we have
\begin{equation}
\label{eqn_RootSpacesRefinedK} \fk_\beta^\ft=\fk_\beta^\ft \cap
\fp \oplus \fk_\beta^\ft\cap \fh
\end{equation}

\begin{equation}
\label{eqn_RootSpacesRefinedM} \fm_\beta^\ft=\fm_\beta^\ft\cap
\fp \oplus \fm_\beta^\ft\cap \fh.
\end{equation}
Furthermore $\fz_\fk(\ft)=\fz_\fk(\ft)\cap \fp\oplus
\fz_\fk(\ft)\cap \fh$ and $\fz_\fm(\ft)=\fz_\fm(\ft)\cap \fp\oplus
\fz_\fm(\ft)\cap \fh$.
\end{Lemma}

\begin{proof}
Let $X\in \fk_\beta^\ft$. According to (\ref{eqn_CartanRefined}),
we can decompose $X$ as $X=X_\fp+X_\fh$, where $X_\fp\in \fk\cap
\fp$ and $X_\fh\in \fk\cap \fh$. First of all we have
\begin{equation*}
-\beta(w)^2X_\fp-\beta(w)^2X_\fh=\ad_w^2(X_\fp)+\ad_w^2(X_\fh).
\end{equation*}
for all $w\in \ft$. Since
\begin{equation*}
\ad_w^2X_\fp\in [\ft,[\ft,\fk\cap \fp]]\subset [\fm\cap
\fp,[\fm\cap \fp,\fk\cap \fp]]\subset [\fm\cap \fp,\fm\cap
\fh]\subset \fk\cap \fp
\end{equation*}
and
\begin{equation*}
\ad_w^2X_\fh\in [\ft,[\ft,\fk\cap \fh]]\subset [\fm\cap
\fp,[\fm\cap \fp,\fk\cap \fh]]\subset [\fm\cap \fp,\fm\cap
\fp]\subset \fk\cap \fh,
\end{equation*}
we conclude $\ad_w^2(X_\fp)=-\beta(w)^2X_\fp$ and
$\ad_w^2(X_\fh)=-\beta(w)^2X_\fh$. Thus we have proven
(\ref{eqn_RootSpacesRefinedK}). The proof of the rest of the
lemma is similar.
\end{proof}

\begin{Rem}
Note that the equations (\ref{eqn_RootSpacesRefinedK}) and
(\ref{eqn_RootSpacesRefinedM}) do not have analogues for the root
spaces $\fk_\alpha$ and $\fm_\alpha$ since they do not
necessarily respect the decomposition (\ref{eqn_CartanRefined}).
\end{Rem}

We can now refine decomposition (\ref{eqn_RootSD_t}) as follows:
\begin{equation}
\label{eqn_RootSpacesRefined_k} \fk=(\fz_\fk(\ft)\cap \fh\oplus
\sum_{\beta\in \Delta^+_\ft}\fk_\beta^\ft\cap
\fh)\oplus(\fz_\fk(\ft)\cap \fp\oplus \sum_{\beta\in
\Delta^+_\ft}\fk_\beta^\ft\cap \fp)
\end{equation}
and
\begin{equation}
\label{eqn_RootSpacesRefined_m} \fm=(\fz_\fm(\ft)\cap \fh\oplus
\sum_{\beta\in \Delta^+_\ft}\fm_\beta^\ft\cap \fh) \oplus
(\ft\oplus \sum_{\beta\in \Delta^+_\ft}\fm_\beta^\ft\cap \fp).
\end{equation}

\section{Eigenvalues of the Shape Operator: Commuting Involutions}
\label{sec_EV_ComInv}

In this section, $H\subset G$ is a symmetric subgroup
corresponding to an involution commuting with the one of $K$.
Recall the refined Cartan decomposition (\ref{eqn_CartanRefined})
\begin{equation*}
\fg=\fk\cap \fp\oplus \fk\cap \fh\oplus \fm\cap \fh \oplus
\fm\cap \fp
\end{equation*}
and that $\fa$ is a maximal abelian subspace of $\fm$ containing
a maximal abelian subspace $\ft$ of $\fm\cap \fp$. Then $\fa$ can
be written as $\fa=\ft\oplus \ft'$, where $\ft'\subset \fm\cap
\fh$.

Let $w\in \ft$ and set $p=\Exp(w)$ (which we do not assume to be
regular) and $\gamma(t)=\Exp(tw)$. Our first goal is to express
the tangent space $T_pHp$ in terms of the restricted roots. Note
that the case $H=K$ in the following proposition is the content
of Proposition 3 of \cite{Verhoczki}.

\begin{Prop}
\label{prop_TangSpaceInRoots} The tangent space $T_pH p$ coincides
with the parallel displacement of $(\fz_\fm(\ft)\cap \fh)\oplus
V_1\oplus V_2\subset \fm$ along $\gamma$, where
\begin{equation*}
V_1=\sum_{\beta\in \Delta^+_\ft,\, \beta(w)\notin
\frac{\pi}{2}+\pi\Z}\fm_\beta^\ft\cap \fh \quad\text{ and }\quad
V_2=\sum_{\beta\in \Delta^+_\ft,\, \beta(w)\notin
\pi\Z}\fm_\beta^\ft\cap \fp.
\end{equation*}
\end{Prop}

\begin{proof}
Regarding the elements of $\fh$ as Killing fields on $M$, the
equations in (\ref{eqn_G_Killing}) yield
\begin{equation}
\label{eqn_H_as_KillingField} \fm\cap \fh=\{X\in \fh\mid (\nabla
X)_{eK}=0\},\quad \fk \cap \fh= \{X\in \fh\mid X(eK)=0\}.
\end{equation}
Of course $T_pHp=U_1+ U_2$, where $U_1= \{X(p)\mid X\in \fh\cap
\fm\}$ and $U_2=\{X(p)\mid X\in \fh\cap \fk\}$.

For any $\alpha\in \Delta^+$, let $\{X^\alpha_i\}_{i\in I_\alpha}$
be an orthonormal basis of $\fm_\alpha$; furthermore let
$\{X^0_i\}_{i\in I_0}$ be an orthonormal basis of $\ft'$. Let
$E^\alpha_i$ and $E^0_i$ be the parallel fields along $\gamma$
with $E^\alpha_i(0)=X^\alpha_i(eK)$ and $E^0_i(0)=X^0_i(eK)$,
respectively.

For $X\in \fh$ let $Y=X|_\gamma$ be the Jacobi field along
$\gamma$ obtained by restricting $X$ to $\gamma$. Since $Y(t)$ is
tangent to the orbit through $\gamma(t)$ for all $t$, it follows
from the description of Jacobi fields on symmetric spaces that

\begin{align}
\label{eqn_JacobiFieldGeneral} Y(t)=&\sum_{i\in I_0}
(a_i+b_it)E^0_i(t)+\nonumber \\
&\sum_{\alpha\in \Delta^+}\sum_{i\in I_\alpha}
(c_i\sin(\alpha(w)t)+d_i\cos(\alpha(w)t))E^\alpha_i(t)
\end{align}
for some constants $a_i,b_i,c_i,d_i$.

Consider first the case $X\in \fm_\beta^\ft\cap \fh$ for some
$\beta\in \Delta^+_\ft$ and let $v=X(eK)$. According to Lemma
\ref{lem_RelationRootSDs}, we can write
\begin{equation*}
X=\sum_{\alpha\in \Delta^+:\, \alpha|_\ft=\pm\beta}\sum_{i\in
I_\alpha} \lambda_{\alpha,i} X^\alpha_i
\end{equation*}
for some constants $\lambda_{\alpha,i}\in \R$. Since $X\in \fm$,
we have $Y'(0)=0$ because of (\ref{eqn_H_as_KillingField}), so we
get $b_i=c_i=0$. It follows that
\begin{align*}
Y(t)&=\sum_{\alpha\in \Delta^+:\, \alpha|_\ft=\pm\beta}\sum_{i\in
I_\alpha}\lambda_{\alpha,i}\cos(\alpha(w)t)E^\alpha_i \\
&=\cos(\beta(w)t)\sum_{\alpha\in \Delta^+:\,
\alpha|_\ft=\pm\beta}\sum_{i\in I_\alpha}\lambda_{\alpha,i}E^\alpha_i(t)\\
&=\cos(\beta(w)t){}_\gamma ||_0^t v.
\end{align*}
We thus have shown that if $\beta(w)\notin \frac{\pi}{2}+\pi\Z$,
then the parallel transport of $\fm_\beta^\ft\cap \fh$ along
$\gamma$ is contained in $U_1$.

Let now $X\in \fz_\fm(\ft)\cap \fh$ and $v=X(eK)$. Lemma
\ref{lem_RelationRootSDs} yields
$$X=\sum_{i\in I_0} \mu_i X^0_i +
\sum_{\alpha\in \Delta^+:\, \alpha|_\ft=0} \lambda_{\alpha,i}
X^\alpha_i$$ for some constants $\mu_i$ and $\lambda_{\alpha,i}$.
We obtain
\begin{align*}
Y(t)=\sum_{i\in I_0}\mu_i E^0_i(t)+\sum_{\alpha\in \Delta^+:\,
\alpha|_\ft=0}\sum_{i\in
I_\alpha}\lambda_{\alpha,i}\cos(\alpha(w)t)E^\alpha_i ={}_\gamma
||_0^t v,
\end{align*}
so the parallel transport of $\fz_\fm(\ft)\cap \fh$ along $\gamma$
is contained in $U_1$. It is now clear that $U_1$ is the direct
sum of the parallel transport of $\fz_\fm(\ft)\cap \fh\oplus V_1$.

It remains to describe $U_2$. For any $\alpha\in \Delta^+$, let
$\{Z_i^\alpha\}_{i\in I_\alpha}$ be the orthonormal basis of
$\fk_\alpha$ which is related to $\{X_i^\alpha\}$ by
$[Z_i^\alpha,u]=\alpha(u)X^\alpha_i$ for all $u\in \fa$. For any
$\alpha\in \Delta^+$, we have
\begin{equation*}
[\fk_\alpha,u]=\begin{cases}\fm_\alpha & \text{if }\alpha(u)\neq
0 \\ 0 & \text{if }\alpha(u)=0.\end{cases}
\end{equation*}
Lemma \ref{lem_RelationRootSDs} now yields that an analogous
relation is true for the root spaces with respect to $\ft$: for
all $\beta\in \Delta^+_\ft$, we have

\begin{equation*}
[\fk^\ft_\beta,u]=\begin{cases}\fm^\ft_\beta & \text{if }\beta(u)\neq 0 \\
0 & \text{if }\beta(u)=0.\end{cases}
\end{equation*}
For $X\in \fk$ we have that
\begin{align*}
Y(t)&=X(\Exp(tw))=\left.\frac{d}{ds}\right|_{s=0}\exp(sX)\exp(tw)K\\
&=\left.\frac{d}{ds}\right|_{s=0}\exp(\Ad_{\exp sX} tw)K = (d\,
\Exp)_{tw}(t[X,w]).
\end{align*}
Hence, those $X\in \fh$ which lie in $\fz_\fk(\ft)\cap \fh$ do
not contribute to $U_2$. We thus have
\begin{equation*}
U_2=(d\,\Exp)_w\sum_{\beta\in \Delta^+_\ft}[\fk_\beta^\ft\cap
\fh,w]=(d\,\Exp)_w\sum_{\beta\in \Delta^+_\ft,\, \beta(w)\neq 0}
\fm_\beta^\ft\cap \fp.
\end{equation*}
Let $v\in \fm_\beta^\ft\cap \fp$, where $\beta(w)\neq 0$, and
write
\begin{equation*}
v=\sum_{\alpha\in \Delta^+:\, \alpha|_\ft=\pm\beta}\sum_{i\in
I_\alpha}\lambda_{\alpha,i} X^\alpha_i.
\end{equation*}
for some constants $\lambda_{\alpha,i}$. Define $X\in
\fk_\beta^\ft\cap \fh$ to be
\begin{equation*}
X=\sum_{\alpha\in \Delta^+:\, \alpha|_\ft=\pm\beta}\sum_{i\in
I_\alpha}\pm\lambda_{\alpha,i} Z^\alpha_i.
\end{equation*}
By definition, we have $[X,w]=\beta(w)v$.

Since $Y$ is the unique Jacobi field along $\gamma$ with $Y(0)=0$
and $Y'(0)=[X,w]=\beta(w)v$, we get
\begin{align*}
Y(t)&=\sum_{\alpha\in \Delta^+:\, \alpha|_\ft=\pm\beta}\sum_{i\in
I_\alpha}\pm\lambda_{\alpha,i}\sin(\alpha(w)t)E^\alpha_i(t)=\sin(\beta(w)t){}_\gamma
||_0^t v.
\end{align*}
It follows that $Y(1)$ vanishes if and only if $\beta(w)\notin
\pi\Z$. We have thus proven that $U_2$ is the parallel
displacement of $V_2$ along $\gamma$.
\end{proof}
\begin{Cor} \label{cor_RegularPointsRoots}
The point $p=\Exp(w)\in \Sigma$ is a regular point of the
$H$-action if and only if
\begin{enumerate}
\item $\beta(w)\notin \frac{\pi}{2}+\pi\Z$ for all $\beta\in \Delta^+_\ft$
with $\fm_\beta^\ft\cap \fh\neq \{0\}$ and
\item $\beta(w)\notin \pi\Z$ for all $\beta\in \Delta^+_\ft$ with
$\fm_\beta\cap \fp\neq \{0\}$.
\end{enumerate}
\end{Cor}
Choose a vector $u\in\ft$, let $c(t)=\Exp(w+tu)$ and
$u(p)=\overset{.}{c} (0)$. By Corollary \ref{Cor_AllCommute}, the
shape and curvature operators can be simultaneously diagonalized.
A concrete such diagonalization is given in the following theorem.
\begin{Thm} \label{thm_ShapeOp}
The decomposition of $T_pHp$ into parallel displacements of the
root spaces described in Proposition \ref{prop_TangSpaceInRoots}
is compatible with the decomposition into the eigenspaces of the
shape operator $A_{u(p)}$ of $Hp$. More precisely,
\begin{enumerate}
\item For $v\in \fm_\beta^\ft\cap \fh$ with $\beta(w)\notin
\frac{\pi}{2}+\pi\Z$, we have
\begin{equation}\label{eqn_ShapeEV_H}
A_{u(p)}({}_\gamma||_0^1v)=\beta(u)\tan(\beta(w)){}_\gamma
||_0^1v.\end{equation}
\item For $v\in \fm_\beta^\ft\cap \fp$ with $\beta(w)\notin \pi\Z$, we
have \begin{equation}\label{eqn_ShapeEV_P}A_{u(p)}({}_\gamma
||_0^1v)=-\beta(u)\cot(\beta(w)){}_\gamma ||_0^1 v.\end{equation}
\item For $v\in \fz_\fm(\ft)\cap \fh$, we have
\begin{equation}\label{eqn_ShapeEV_T}A_{u(p)}({}_\gamma
||_0^1v)=0.\end{equation}
\end{enumerate}
\end{Thm}
\begin{proof}
For any $s\in [0,1]$, let $\gamma_s(t):=\Exp(t(w+su))$. Note that
$\gamma_0=\gamma$.

First of all, let $X\in \fm_\beta^\ft\cap \fh$, where
$\beta(w)\notin \frac{\pi}{2}+\pi\Z$, and set $v=X(eK)$. Let
$Y_s$ be the Jacobi field obtained by restriction of the Killing
field $X$ to $\gamma_s$. The initial values of $Y_s$ are
$Y_s(0)=v$ and $Y_s'(0)=0$; as in the proof of Proposition
\ref{prop_TangSpaceInRoots}, we get
\begin{equation}
Y_s(t)=\cos(\beta(w+su)t){}_{\gamma_s} ||_0^1 v,
\end{equation}
since $v$ is contained in a sum of root spaces corresponding to
roots whose restrictions to $\ft$ coincide. We are interested in
the $Hp$-Jacobi field $Y(t)=Y_t(1)$ along $c$. Its initial values
are $Y(0)=\cos(\beta(w)){}_{\gamma} ||_0^1 v$ and
\begin{align}
Y'(0)&=\left.\frac{d}{dt}\right|_{t=0}\cos(\beta(w+tu)){}_{\gamma_t}
||_0^1 v \nonumber \\
&= -\beta(u)\sin(\beta(w)){}_{\gamma} ||_0^1 v,
\end{align}
since $\left.\frac{\nabla}{dt}\right|_{t=0}{}_{\gamma_t}
||_0^1v=0$ (use Lemma 8.3.2 of \cite{BCO}, together with the fact
that the $\gamma_s$ lie in the flat section $\Sigma$). The fact
that $Y$ is an $Hp$-Jacobi field along $c$ now implies
\begin{align*}
Y'(0)+&A_{u(p)}Y(0)\\
&=-\beta(u)\sin(\beta(w)){}_\gamma
||_0^1v+\cos(\beta(w))A_{u(p)}({}_\gamma ||_0^1v)\in \nu_pHp,
\end{align*} so
we get \begin{equation} A_{u(p)}({}_\gamma ||_0^1
v)=\beta(u)\tan(\beta(w)){}_\gamma ||_0^1 v,\end{equation} which
is Equation (\ref{eqn_ShapeEV_H}).

In order to prove Equation (\ref{eqn_ShapeEV_P}), choose some
vector $v\in \fm_\beta^\ft\cap \fp$ with $\beta(w)\notin \pi\Z$.
Let $X\in \fk_\beta^\ft$ be such that $[H,X]=-\beta(H)v$ for all
$H\in \ft$. We have $X\in \fh$ since $\beta(w)\neq 0$. Now
continue exactly as above: let $Y_s$ be the Jacobi field obtained
by restriction of the Killing field $X$ to $\gamma_s$. Its initial
values are $Y_s(0)=0$ and $Y_s'(0)=[X,w+su]=\beta(w+su)v$, so we
get

\begin{equation}
Y_s(t)=\sin(\beta(w+su)t){}_{\gamma_s} ||_0^tv.
\end{equation}
The $Hp$-Jacobi field $Y(t)=Y_t(1)$ has initial values
$Y(0)=\sin(\beta(w)){}_\gamma ||_0^1v$ and
\begin{align}
Y'(0)&=\left.\frac{d}{dt}\right|_{t=0}\sin(\beta(w+tu)){}_{\gamma_t}
||_0^1v \nonumber\\
&=\alpha(u)\cos(\beta(w)){}_\gamma ||_0^1v, \end{align} so we
obtain
\begin{align*}
Y'(0)+&A_{u(p)}Y(0)\\
&=\beta(u)\cos(\beta(w)){}_\gamma
||_0^1v+\sin(\beta(w))A_{u(p)}({}_\gamma ||_0^1v)\in \nu_pHp;
\end{align*}
thus, Equation (\ref{eqn_ShapeEV_P}) follows.

Finally, let $X\in \fz_\fm(\ft)\cap \fh'$ and set $v=X(eK)$. Let
again $Y_s$ be the restriction of $X$ along $\gamma_s$ and
$Y(t)=Y_t(1)$. Then we see that $Y_s(t)={}_{\gamma_s} ||_0^1v$ and
hence $Y$ satisfies the initial conditions $Y(0)={}_\gamma ||_0^1
v$ and $Y'(0)=0$. Equation (\ref{eqn_ShapeEV_T}) follows
immediately.

\end{proof}

\section{Eigenvalues of the Shape Operator: General Case}
\label{sec_EV_General}

In this section we will determine, as far as possible, the
eigenvalues of the shape operators in the general case of an
arbitrary Hermann action. This is independent of the calculations
in section \ref{sec_EV_ComInv}.

The following proposition shows the dependence of the normal
direction. Let the origin be a regular point, denoted by $p$.
Note that $T_pHp\cong \pr_\fm \fh=\sum_\beta \fm_\beta^\ft \oplus
(\fz_\fm(\ft)\cap \pr_\fm \fh)$.
\begin{Prop}\label{prop_ShapeOpGeneral}
\text{\emph{(1)}} There exists a refinement $\fm_\beta^\ft=\sum_i
V_{\beta,i}$ of the root spaces, together with constants
$c_{\beta,i}$, such that for $v\in \nu_pHp$ and all $x\in
V_{\beta,i}$, we have
$$A_vx=c_{\beta,i}\beta(v)x.$$
\text{\quad\emph{(2)}} For $x\in \fz_\fm(\ft)\cap \pr_\fm\fh$, we
have
$$A_vx=0.$$
\end{Prop}
\begin{proof}
Let $x\in T_pHp$ be a common eigenvector of all curvature and
shape operators in normal directions.

Let the linear form $f:\nu_pHp\to \R$, depending on $x$, be
defined by $A_vx=f(v)x$. Choose $X\in \fh$ such that $X(p)=x$. We
write $X=X_\fk+X_\fm$ with $X_\fk\in \fk$ and $X_\fm\in \fm$.

For $v\in\nu_pHp$, we denote by $\xi'\in \fm$ the Killing vector
field induced by transvections with $\xi'(p)=v$ and by $\xi$ the
$H$-equivariant parallel normal vector field with $\xi(p)=v$.
Then we have
\begin{align*}
[\xi',X_\fk](p)&=-\nabla_vX_\fk+\nabla_{\underbrace{X_\fk(p)}_{=0}}\xi'=-\nabla_vX=-\nabla_x\xi-\underbrace{[\xi,X]}_{=0}\\
&=A_vx=f(v)x;
\end{align*}
since $[\xi',X_\fk]\in [\fm,\fk]\subset \fm$, it follows that
\begin{equation}[\xi',X_\fk]=f(v)X_\fm.\label{Eqn_XkXpBracket}\end{equation}

Let us now first regard the case of $x\in \fm_\beta^\ft$ for some
root $\beta$. Since (\ref{Eqn_XkXpBracket}) is valid for all
$\xi'\in \fm$, we can write $X_\fk=X_{\fk,0}+X_{\fk,1}$ with
$X_{\fk,0}\in \fz_\fk(\fm)$ and $X_{\fk,1}\in \fk_\beta^\ft$; the
vector $Z\in \fm_\beta^\ft$ related to $X_{\fk,1}$ is a multiple
of $X_\fm$, i.e.~there exists some constant $c$, independent of
$v$, with $f(v)=c\cdot \alpha(v)$.

If $x\in \fm$ is such that $[X_\fm,\ft]=0$, it follows from
(\ref{Eqn_XkXpBracket}) that
$$f(v)\left<X_\fm,X_\fm\right>=\left<[\xi',X_\fk],X_\fm\right>=-\left<X_\fk,[\xi',X_\fm]\right>=0;$$
hence $A_vx=0.$
\end{proof}
\begin{Rem} The explicit description of the eigenvalues in the case of
commuting involutions was possible because there existed a point
$p\in M$ such that every $H$-Killing vector field could be
written as the sum of one vanishing at $p$ and one with
derivative vanishing at $p$, leading to eigenvalues containing
either a cotangent or a tangent.

In the general case such a point does not exist, but for each
Killing vector field we can choose a point where either it
vanishes or its derivative. For the following calculation we
choose to express the $c_{\beta,i}$ in terms of zeros of the
Killing vector fields themselves; hence, only the cotangent
occurs.
\end{Rem}

Now we can investigate how the $c_{\beta,i}$ and the eigenvalues
of the shape operators change when varying the orbit. In the
notation above, let $Y_w:=X\circ\gamma_w$, where
$\gamma_w(t)=\Exp(tw)$ is the geodesic in direction $w$. For
$x=X(p)\in V_{\beta,i}\subset \fm_\beta^\ft$, we have
$$Y_w(t)=(-c_{\beta,i} \sin(\beta(w)t)+\cos(\beta(w)t)){}_{\gamma_w}||_0^t
x.$$ If $t_{\beta,i}$ is a fixed zero of some $Y_{w_0}$ with
$\beta(w_0)=1$, we can write $c_{\beta,i}=\cot(t_{\beta,i})$.
Regarding the $H\gamma_w(1)$-Jacobi field $Y(t):=Y_{w+tv}(1)$, we
can determine the shape operator of the orbit $H\gamma_w(1)$:
\begin{align*}
0&=Y'(0)+A_vY(0)\\
&=(-c_{\beta,i}\beta(v)\cos(\beta(w))-\beta(v)\sin(\beta(w))){}_{\gamma_w}||_0^1
x+A_v Y_w(1).
\end{align*}
Hence
\begin{align*}
A_v({}_{\gamma_w}||_0^1
x)&=\beta(v)\frac{c_{\beta,i}\cos(\beta(w))+\sin(\beta(w))}{-c_{\beta,i}\sin(\beta(w))+\cos(\beta(w))}{}_{\gamma_w}||_0^1x\\
&=\beta(v)\frac{\cos(t_{\beta,i})\cos(\beta(w))+\sin(t_{\beta,i})\sin(\beta(w))}{-\cos(t_{\beta,i})\sin(\beta(w))+\sin(t_{\beta,i})\cos(\beta(w))}{}_{\gamma_w}||_0^1x\\
&=\beta(v)\cot(t_{\beta,i}-\beta(w)){}_{\gamma_w}||_0^1x.
\end{align*}

\section{Applications}
\label{sec_Appl}
 Let $H\subset G$ be a symmetric subgroup such that the corresponding
involutions commute. Let $p\in M$ be regular and
$\Sigma=\Exp(\ft)$ be the section through $p$. Denote the
generalized Weyl group of the action by $W$ and define a function
$\vartheta:\ft\to \R$ by

\begin{equation}
\vartheta(w)=\prod_{\beta\in
\Delta^+_\ft}|\sin(\beta(w))|^{p_\beta}|\cos(\beta(w))|^{h_\beta},
\end{equation}
where the exponents $p_\beta$ and $h_\beta$ are the relative root
multiplicities defined by $p_\beta=\dim \fm_\beta^\ft\cap \fp$ and
$h_\beta=\dim \fm_\beta^\ft\cap \fh$. Since $\vartheta$ is
invariant under the reflections in the singular hyperplanes in
$\ft$, it can be regarded as a function on $\Sigma$. We will
reprove the following theorem from \cite{F-J} using our
calculation of the shape operators. The proof is similar to the
arguments in \cite{CNV}, where the case $H=K$ is treated.

\begin{Thm}\label{thm_IntFormula}
For any integrable function $f$ on $M$, we have
\begin{equation}
\label{eqn_IntFormulaNotInvariant} \int_M
f(x)\,dx=\frac{1}{|W|\cdot\vartheta(p)}\int_\Sigma\left(\int_{H/H_p}f(hq)\,
d(hH_p)\right)\vartheta(q)\, dq,
\end{equation}
where the Riemannian measure on $H/H_p$ is chosen to be the one
induced by $H\cdot p\subset M$.

If $f$ is additionally $H$-invariant, we have

\begin{equation}
\label{eqn_IntFormulaInvariant} \int_M
f(x)\,dx=\frac{\Vol(M)}{\int_\Sigma \vartheta(q)\, dq}
\int_\Sigma f(q)\vartheta(q)\, dq.\end{equation}
\end{Thm}

\begin{Rem}
If $H$ is an arbitrary symmetric subgroup of $G$, the theorem is
true with $\vartheta$ defined as follows: Conjugate $H$ such that
the regular point $p$ is the origin $eK\in M$ and consider the
decomposition
$$T_pHp=(\fz_\fm(\fh)\cap \pr_\fm \fh)\oplus\sum_{\beta\in
\Delta^+_\ft}\sum_i V_{\beta,i}$$ of Proposition
\ref{prop_ShapeOpGeneral}. Choose nonzero $x_{\beta,i}\in
V_{\beta,i}$ and $X_{\beta,i}\in \fh$ with
$X_{\beta,i}(p)=x_{\beta,i}$. Then choose $w_{\beta,i}\in \ft$
with $\beta(w_{\beta,i})=1$ and let $t_{\beta,i}$ be a zero of
the Jacobi field $X_{\beta,i}\circ \gamma_{\beta,i}$, where
$\gamma_{\beta,i}$ is the geodesic in direction $w_{\beta,i}$.
Then define $\vartheta:\ft\to \R$ as
\begin{equation}
\vartheta(w)=\prod_{\beta\in \Delta_\ft^+}\prod_i
|\sin(\beta(w)-t_{\beta,i})|^{\dim V_{\beta,i}}.
\end{equation}
Since the proof is completely analogous using the calculations in
Section \ref{sec_EV_General} instead of Theorem
\ref{thm_ShapeOp}, we will prove the theorem only in the case of
commuting involutions.

We also remark that in the noncompact case, the theorem remains
true if $\vartheta$ is defined using hyperbolic functions.
\end{Rem}

\begin{Lemma} \label{lem_Fubini} {\bf (Generalized Cavalieri Principle)}
Let $M$ be a Riemannian manifold such that a subset $U\subset M$
of full measure can be written as $L\times N$, equipped with a
Riemannian metric of the form $g(q,r)=\left(\begin{matrix}h(q,r) &
0 \\ 0 & k(r)\end{matrix}\right)$, where $(q,r)\in L\times N$.
Then for any integrable function $f$ on $M$, we have
\begin{equation*} \int_M f(p)\, dM=\int_{L\times N} f(q,r)\,
d(L\times N)=\int_N\left( \int_{L_r} f(q,r)\, dL_r\right) dN.
\end{equation*}
\end{Lemma}
\begin{proof}
Applying Fubini's theorem in coordinates $\phi:V\times
W\overset{\sim}{\to} V'\times W'\subset U$ yields
\begin{align*}
\int_{V'\times W'}&f(q,r)\, dM=\int_{V\times W} f\circ
\phi(x,y)\sqrt{\det
(g_{ij}(x,y))}\, d(x,y)\\
&=\int_{W}\left(\int_V f\circ \phi(x,y)\sqrt{\det (h_{ij}(x,y))}\,
dx\right)\sqrt{\det (k_{ij}(y))}\, dy.\\
&= \int_{W'}\left(\int_{V'_r}f(q,r)\, dV'_r\right)dW'.
\end{align*}\vskip-20pt
\end{proof}

Applying Lemma \ref{lem_Fubini} to the regular set in $M$, we
obtain that for any integrable function $f$ on $M$,
\begin{align}
\int_M f(x)\, dx&=\int_{\bar{Q}} \left(\int_{H\cdot \Exp(w)}
f(x)\,
dx\right) \, dw \nonumber \\
&=\frac{1}{|W|}\int_\Sigma\left( \int_{Hq}f(x)\, dx\right)\, dq,
\label{eqn_IntSigmaOrbit}
\end{align}
where $\bar{Q}\subset \ft$ is a generalized Weyl chamber, and $W$
is the generalized Weyl group of the $H$-action on $M$.

\begin{Lemma} \label{lem_CompareVolumes}
The Riemannian densities $\mu$ and $\mu'$ of the $H$-orbits
through two regular points $p=\Exp(w)$ and $q=\Exp(w+tu)$ are
related by the formula
\begin{align*}
\mu'\circ \varphi=F_{p}(q)\cdot \mu,
\end{align*} where
\begin{align}
F_{p}(q)= &\prod_{\beta\in \Delta^+_\ft}
\left(\cos(\beta(tu))+\cot(\beta(w))\sin(\beta(tu))\right)^{p_\beta}
\nonumber \\ & \quad\,\times
\left(\cos(\beta(tu))-\tan(\beta(w))\sin(\beta(tu))\right)^{h_\beta}.\label{eqn_RelativeDensities}
\end{align}
Here the exponents are the relative root multiplicities
$h_\beta=\dim \fm_\beta^\ft\cap \fh$ and $p_\beta=\dim
\fm_\beta^\ft\cap \fp$.
\end{Lemma}
\begin{proof}
Let $\xi$ be the $H$-equivariant vector field on $Hp$ determined
by $\xi(p)=u(p)$. The polarity of the $H$-action implies that
$\xi$ is a parallel normal vector field (\cite{PalaisTerng},
Theorem 5.7.1) and the map $\varphi:=\Exp_{Hp}\circ (t\xi)$ is a
diffeomorphism between $Hp$ and $Hq$. Since $\xi$ is
$H$-equivariant, so is $\varphi$, and consequently, for $X\in
\fh$, we have $d\varphi(X\cdot p)=X\cdot q$. In other words,
applying $d\varphi$ to $X\cdot p$ is the same as evaluating the
Jacobi field obtained by restriction of the Killing vector field
$X$ to the geodesic $c(t)=\Exp_{Hp}(tu(p))$ at time $t$.

In order to compare the Riemannian densities of $Hp$ and $Hq$, we
need to calculate the determinant of $d\varphi_p$; the
argumentation above shows that this amounts to explicitly
calculate these Jacobi fields.

Let $X\in \fm_\beta^\ft\cap \fh$, $v= X(eK)$ and $v'={}_\gamma
||_0^1v\in T_pHp$. The initial values of the Jacobi field
$Y_v(t)=X\cdot c(t)$ along $c$ have been calculated in the proof
of Theorem \ref{thm_ShapeOp}. Therefore
\begin{equation}\label{eqn_NormalizedJacFieldH}
\frac{1}{\cos(\beta(w))}Y_v(t)=(\cos(\beta(tu))-\tan(\beta(w))\sin(\beta(tu))){}_c||_0^t
v'.\end{equation}

If $X\in \fz_\fm(\ft)\cap \fh$, $v=X(eK)$ and $v'={}_\gamma ||_0^1
v$, then \begin{equation}
\label{eqn_NormalizedJacFieldT}Y_v(t)={}_c ||_0^t v'\end{equation}
is the restriction of $X$ to $c$.

If $v\in \fm_\beta^\ft\cap \fp$, $v'={}_\gamma ||_0^1v$ and $X\in
\fk_\beta^\ft\cap \fh$ is such that $[H,X]=-\beta(H)v$ for all
$H\in \ft$, we analogously see that the Jacobi field $Y_v$,
defined by $Y_v(t)=X\cdot c(t)$, is given by
\begin{equation}
\label{eqn_NormalizedJacFieldP}
\frac{1}{\sin(\beta(w))}Y_v(t)=(\cos(\beta(tu))+\cot(\beta(w))\sin(\beta(tu))){}_c||_0^t
v'.\end{equation}

If $\mu$ and $\mu'$ are the Riemannian densities of $Hp$ and $Hq$,
respectively, we conclude from (\ref{eqn_NormalizedJacFieldH}),
(\ref{eqn_NormalizedJacFieldT}) and
(\ref{eqn_NormalizedJacFieldP}) that
\begin{align*}\mu'(\varphi(p))=\mu(p)\cdot &\prod_{\beta\in
\Delta^+_\ft}
\left(\cos(\beta(tu))+\cot(\beta(w))\sin(\beta(tu))\right)^{p_\beta}
\\ &\quad\,\times
\left(\cos(\beta(tu))-\tan(\beta(w))\sin(\beta(tu))\right)^{h_\beta},\end{align*}
since the Jacobi fields in question are multiples of parallel
vector fields and therefore stay orthogonal to each other.
\end{proof}
\begin{Cor}
\label{cor_CompareVolumes} With the same notation as in Lemma
\ref{lem_CompareVolumes}, the volumes of $Hp$ and $Hq$ are
related by \begin{equation*} \Vol(Hq)=\Vol(Hp)\cdot
F_p(q).\end{equation*}
\end{Cor}

Lemma \ref{lem_CompareVolumes} now enables us to replace the inner
integral in Equation (\ref{eqn_IntSigmaOrbit}) by an integral over
a fixed regular orbit $Hp$, where $p\in\Sigma$. If we denote the
Riemannian density of the orbit $Hq$ by $\mu_q$, and
$\varphi_q:Hp\to Hq$ is the equivariant diffeomorphism introduced
in Lemma \ref{lem_CompareVolumes}, we have
\begin{align}
\int_{M} f(x)\,
dx&=\frac{1}{|W|}\int_{\Sigma}\left(\int_{Hq}f(x)\,
dx\right)\, dq \nonumber\\
&=\frac{1}{|W|}\int_{\Sigma}
\left(\int_{Hp}f(\varphi_q(x))\frac{\mu_q(\varphi_q(x))}{\mu_p(x)}\,
dx\right) \, dq\nonumber \\
&=\frac{1}{|W|}\int_{\Sigma} \left(\int_{H/H_p}f(hq)\,
dh\right)F_{p}(q)\, dq.\label{eqn_IntOverHp}
\end{align}

A short calculation using the addition theorems for the
trigonometric functions shows that
\begin{equation}\label{eqn_DensitiesTheta}
\frac{\mu_q\circ
\varphi_q}{\mu_p}=F_p(q)=\frac{\vartheta(q)}{\vartheta(p)}.
\end{equation}
Now equality (\ref{eqn_IntFormulaNotInvariant}) follows from
(\ref{eqn_IntOverHp}) and (\ref{eqn_DensitiesTheta}).

If we additionally assume $f$ to be $H$-invariant, we can argue
as in Theorem 3.5 of \cite{CNV}: on the one hand, we clearly have
\begin{equation*}
\int_M f(x)\, dx= \frac{1}{|W|}\int_\Sigma \Vol(Hq)f(q)\, dq.
\end{equation*}
On the other hand, Corollary \ref{cor_CompareVolumes} and Equation
(\ref{eqn_DensitiesTheta}) show that
\begin{equation*}
\Vol(Hq)=\frac{\Vol(Hp)}{\vartheta(p)}\cdot \vartheta(q),
\end{equation*}
so the quotient $\Vol(Hp)/\vartheta(p)$ does not depend on $p$,
and we may denote it by $V(M)$. Then we have
\begin{equation*}
\Vol(M)=\frac{1}{|W|}\int_\Sigma \Vol(Hq)\,
dq=\frac{V(M)}{|W|}\cdot \int_\Sigma \vartheta(q)\, dq.
\end{equation*}

Combining these three equalities, we have shown
(\ref{eqn_IntFormulaInvariant}). On the way, we have also proven
\begin{Prop} For any regular point $p$,
\begin{equation}
\Vol(Hp)=\Vol(M)|W|\frac{\vartheta(p)}{\int_\Sigma
\vartheta(q)dq}.
\end{equation}
\end{Prop}

\section{Example: $\U(p+q)$ Acting on $\SO(2p+2q)/\rmS(\rmO(2p)\times
\rmO(2q))$} \label{sec_Example}

In this section, we determine the function $\vartheta$ introduced
in Section \ref{sec_Appl} for the $H=\U(p+q)$-action on the
symmetric space
$$M=\SO(2p+2q)/\rmS(\rmO(2p)\times \rmO(2q)).$$ Let us assume that
$p\le q$.

On the Lie algebra $\fg=\fso(2p+2q)$, we consider the two
involutions $\sigma_1(X)=I_{2p,2q}XI_{2p,2q}$ and
$\sigma_2(X)=J_{p+q}XJ_{p+q}^{-1}$, where
$$\quad I_{2p,2q}=\left(\begin{matrix}I_p & 0 & 0 & 0 \\ 0 & I_p & 0 & 0
\\ 0 & 0 & -I_q & 0 \\ 0 & 0 & 0 & -I_q\end{matrix}\right) \quad \text{
and }    J_{p+q}= \left(\begin{matrix}0 & I_p & 0 & 0 \\ -I_p & 0
& 0 & 0
\\ 0 & 0 & 0 & I_q \\ 0 & 0 & -I_q & 0\end{matrix}\right). $$ Clearly,
$\sigma_1$ and $\sigma_2$ commute. The corresponding fixed point
algebras are $\fk=\fg^{\sigma_1}=\fso(2p)\oplus \fso(2q)$ and

\begin{equation*}
\fh=\fg^{\sigma_2}=\left\{\left.\left( \begin{array}{cc|cc} a & -b
& c & -d \\ b & a & d & c \\ \hline -c^t & -d^t & e & -f \\
\smash{\underbrace{d^t}_{p}} & \smash{\underbrace{-c^t}_{p}} &
\smash{\underbrace{f}_{q}} & \smash{\underbrace{e}_{q}}
\end{array}\right)\right| \begin{matrix} a^t=-a & e^t=-e \\ b^t=b
& f^t=f \end{matrix} \right\}\cong\fu(p+q).
\end{equation*}\vskip12pt

Let the respective Cartan decompositions of $\fg$ be
$\fg=\fk\oplus \fm$ and $\fg=\fh\oplus \fp$. We have
$$\fm=\left\{\left(\begin{array}{c|c} 0 & a \\ \hline -a^t &
0\end{array}\right)\in \fg\right\}$$ and
\begin{equation*}
\fp=\left\{\left.\left(\begin{array}{cc|cc} u & v & w & x \\ v &
-u & x & -w \\ \hline -w^t & -x^t & y & z \\ -x^t & w^t & z & -y
\end{array}\right)\right| \begin{matrix}u^t=-u & v^t=-v \\ y^t=-y
& z^t=-z\end{matrix}\right\}.
\end{equation*}
The intersections $\fm \cap \fp$ and $\fm \cap \fh$ are given by
$$\fm\cap \fp=\left\{\left(\begin{array}{cc|cc}0 & 0 & w & x \\ 0
& 0 & x & -w \\ \hline -w^t & -x^t & 0 & 0 \\ -x^t & w^t &  0 &
0\end{array}\right)\in \fg\right\}$$ and $$\fm\cap
\fh=\left\{\left(\begin{array}{cc|cc} 0 & 0 & c & -d \\ 0 & 0 & d
& c \\ \hline -c^t & -d^t & 0 & 0 \\ d^t & -c^t &  0 &
0\end{array}\right)\in \fg\right\}.$$ For $1\le i\le p$ and $1\le
j\le q$, let $E_{i,j}$ be the $p\times q$-matrix having $1$ in the
$(i,j)$-th entry and zeros elsewhere, and define

$$Q_{i,j}=\left(\begin{array}{cc|cc} 0 & 0 & E_{i,j} & 0 \\ 0 & 0
& 0 & -E_{i,j}\\ \hline -E_{i,j}^t & 0 & 0 & 0 \\ 0 & E_{i,j}^t &
0 & 0\end{array}\right)\in \fm\cap \fp,$$
$$
R_{i,j}=\left(\begin{array}{cc|cc} 0 & 0 & 0 & E_{i,j} \\ 0 & 0 &
E_{i,j} & 0\\ \hline 0 & -E_{i,j}^t & 0 & 0 \\ -E_{i,j}^t & 0 & 0
& 0\end{array}\right)\in \fm\cap \fp,$$

$$F_{i,j}=\left(\begin{array}{cc|cc} 0 & 0 & E_{i,j} & 0 \\ 0 & 0 &
0 & E_{i,j}\\ \hline -E_{i,j}^t & 0 & 0 & 0 \\ 0 & -E_{i,j}^t & 0
& 0\end{array}\right)\in \fm\cap \fh,$$
$$
G_{i,j}=\left(\begin{array}{cc|cc} 0 & 0 & 0 & -E_{i,j} \\ 0 & 0
& E_{i,j} & 0\\ \hline 0 & -E_{i,j}^t & 0 & 0 \\ E_{i,j}^t & 0 &
0 & 0\end{array}\right)\in \fm\cap \fh.$$

Then we define abelian subalgebras by $$\ft=\sum_{i=1}^p \R
Q_{i,i}\subset \fm\cap \fp,\quad \ft'=\sum_{i=1}^p \R
F_{i,i}\subset \fm\cap \fh,\quad \fa=\ft\oplus \ft'\subset \fm$$
and consider the root space decomposition of $(\fg,\fk)$ with
respect to $\fa$: For $1\le i\le 2p$, let linear forms
$\lambda_i:\fa\to \R$ be defined by
$$
\lambda_i(Q_{j,j})=\begin{cases}\delta_{i,j} & \text{if }i\le p \\
-\delta_{i-p,j} & \text{if }i>p \end{cases},\quad
\lambda_i(F_{j,j})=\begin{cases}\delta_{i,j} & \text{if }i\le p \\
\delta_{i-p,j} & \text{if }i>p \end{cases}.
$$
The restricted roots are $\pm \lambda_i$ for $1\le i\le 2p$ and
$\pm (\lambda_i\pm \lambda_j)$ for $ 1\le i<j\le 2p$; we choose
as positive roots the $\lambda_i$ and the $\lambda_i\pm
\lambda_j$ for $i<j$. Denoting the restriction of $\lambda_i$ to
$\ft$ by $\lambda_i^\ft$, we obtain the adapted root space
decomposition with respect to $\ft$
\begin{equation*}
\fm=\fz_\fm(\ft)\oplus \sum_{1\le i\le
p}\fm_{\lambda_i^\ft}^\ft\oplus \underset{i<j\le 2p}{\sum_{1\le
i\le p}} \fm_{\lambda_i^\ft\pm \lambda_j^\ft}^\ft
\end{equation*}
where
\begin{align*}
\fm_{\lambda_i^\ft}^\ft= \underbrace{\sum_{k=p+1}^q \R
F_{i,k}\oplus \R G_{i,k}}_{\subset \fm\cap\fh} \oplus
\underbrace{\sum_{k=p+1}^q \R Q_{i,k}\oplus \R R_{i,k}}_{\subset
\fm\cap \fp}
\end{align*}
for $1\le i\le p$,
\begin{equation*}
\fm_{\lambda_i^\ft\pm \lambda_j^\ft}^\ft= \underbrace{\R
(F_{i,j}\pm F_{j,i})}_{\subset \fm\cap \fh} \oplus \underbrace{\R
(Q_{i,j}\pm Q_{j,i})}_{\subset \fm\cap \fp}
\end{equation*}
for $1\le i<j\le p$ and
\begin{equation*}
\fm_{\lambda_i^\ft\pm \lambda_j}^\ft=\underbrace{\R (G_{i,j-p}\pm
G_{j-p,i})}_{\subset \fm\cap \fh} \oplus \underbrace{\R
(R_{i,j-p}\pm R_{j-p,i})}_{\subset \fm\cap \fp}
\end{equation*}
for $1\le i\le p$ and $p< j\le 2p$. Therefore,
\begin{align*}
\vartheta&=\prod_{i=1}^p\prod_{j=1}^{2p}
|\sin(\lambda_i^\ft+\lambda_j^\ft)||\sin(\lambda_i^\ft-\lambda_j^\ft)||\cos(\lambda_i^\ft+\lambda_j^\ft)||\cos(\lambda_i^\ft-\lambda_j^\ft)|
\\
& \quad\times \prod_{i=1}^p |\sin (\lambda_i^\ft)|^{2(q-p)} |\cos
(\lambda_i^\ft)
|^{2(q-p)} \\
& = \prod_{i,j=1}^p
|\sin(\lambda_i^\ft+\lambda_j^\ft)|^2|\sin(\lambda_i^\ft-\lambda_j^\ft)|^2|\cos(\lambda_i^\ft+\lambda_j^\ft)|^2|\cos(\lambda_i^\ft-\lambda_j^\ft)|^2
\\
& \quad\times \prod_{i=1}^p |\sin (\lambda_i^\ft)|^{2(q-p)} |\cos
(\lambda_i^\ft) |^{2(q-p)},
\end{align*}
where we have used that $\lambda_{i+p}^\ft=-\lambda_{i}^\ft$ for
$1\le i\le p$.

\end{document}